\documentclass[
submission
]{dmtcs-episciences}

\usepackage[utf8]{inputenc}
\usepackage{subfigure}

\usepackage{mathrsfs}
\usepackage{amsfonts}

\usepackage{amsfonts, amsmath, amssymb}
\usepackage{amssymb,amsfonts,amsmath,
latexsym, epsfig,cite, psfrag,eepic,colordvi,tikz}
\usepackage{amscd,graphics}

\newtheorem{theorem}{Theorem}[section]
\newtheorem{lemma}[theorem]{Lemma}

\usepackage[round]{natbib}

\author{Hongliang Lu\affiliationmark{1}\thanks{Supported  by the National Natural
Science Foundation of China under grant No.11871391 and
Fundamental Research Funds for the Central Universities.}
  \and Wei Wang\affiliationmark{1}
  \and Juan Yan\affiliationmark{2} \thanks{Supported by the National Natural
Science Foundation of China under grant No.11801487, 11971274 and QD1919}}
\title[Anti-factors of  Regular Bipartite Graphs]{Anti-factors of  Regular Bipartite Graphs}
\affiliation{
  School of Mathematics and Statistics, Xi'an Jiaotong University, Xi'an, Shaanxi, 710049, China\\
  Department of Mathematics, Lishui University, Lishui, Zhejiang, 323000, China}

\received{2017-3-30}
\revised{2018-5-10,2018-11-19,2019-4-4,2019-10-10,2020-3-26,2020-5-9}
\accepted{2020-5-9}

\begin{document}
\publicationdetails{22}{2020}{1}{16}{3233}
\maketitle
\begin{abstract}
      Let $G=(X,Y;E)$ be a bipartite graph, where $X$ and $Y$ are color classes and $E$ is the set of edges of $G$.
Lov\'asz and Plummer 
 asked whether one can decide in polynomial time that a given bipartite graph $G=(X,Y; E)$ admits a 1-anti-factor, that is subset $F$ of $E$ such that $d_F(v)=1$ for all $v\in X$ and $d_F(v)\neq 1$ for all $v\in Y$. Cornu\'ejols 
 answered this question in the affirmative. Yu and Liu 
 asked whether, for a given integer $k\geq 3$, every $k$-regular bipartite graph contains a 1-anti-factor. This paper answers this question in the affirmative.
\end{abstract}

\section{Introduction}
\label{sec:in}

In this paper, we consider finite undirected graphs without loops
and multiple edges.
Let $G=(V(G),E(G))$ be a graph with vertex set $V(G)$ and edge set $E(G)$. A graph $G'$ is called a \emph{spanning subgraph} of $G$ if $V(G)=V(G')$ and $E(G')\subseteq E(G)$.
The degree of a vertex  $x$ in
$G$ is denoted by $d_G(x)$, and the set of vertices adjacent to $x$
in $G$ is denoted by $N_G(x)$.  For $x\in V(G)$, we write $N_G[x]=N_G(x)\cup \{x\}$. For $xy\notin E(G)$, $G+xy$ denotes the graph with vertex set $V(G)\cup \{x,y\}$ and edge set $E(G)\cup \{xy\}$.
 For $S\subseteq V(G)$, the subgraph
of $G$ induced by $S$ is denoted by $G[S]$ and $G-S= G[V(G)-S]$. For two disjoint subsets $S,T\subseteq V(G)$, let $E_G(S, T)$ denote the set of edges of $G$ joining
$S$ to $T$ and let $e_G(S,T)=|E_G(S,T)|$. For a positive integer $r$, let $[r]=\{0,1,\ldots,r\}$.  Let $c(G)$ denote the number of connected components of $G$.

Let $G$ be a graph, and for every vertex $x\in V(G)$, let $H(x)$ be a set of
integers. An
\emph{$H$-factor} is a spanning graph $F$ such that
\begin{align}\label{cond_fac}
 d_{F}(x)\in
H(x)\quad \mbox{for all $x\in V(G)$}.
\end{align}
A \emph{matching}  of a graph  is a set of edges   such that no two edges share a vertex in common. A \emph{perfect matching} of a graph is a matching covering all vertices. Clearly, a  matching (or perfect matching) of a graph is also a $\{0,1\}$-factor (1-factor, respectively).
On $1$-factors of bipartite graphs, Hall obtained the following result.

\begin{theorem}[\cite{Hall35}]\label{Hall} Let $k\geq 1$ be an integer.
Every $k$-regular bipartite graph contains a 1-factor.
\end{theorem}

 A spanning subgraph $F$ of bipartite graph $G=(X,Y;E)$  is called a \emph{1-anti-factor} if $d_F(x)=1$ for all $x\in X$ and $d_F(y)\neq 1$ for all $y\in Y$. 
Lov\'asz and Plummer (see \cite{LoPl86}, Page 390) proposed the following problem: can one decide in polynomial time whether a given bipartite graph  admits a 1-anti-factor?

 A set $\{h_1,h_2,\ldots,h_m\}$ of increasing integers is called {\em allowed} (see~\cite{Lovasz})
if $h_{i+1}-h_i\le2$ for all $1\le i\le m-1$.
Let $H\colon V(G)\to 2^{\mathbb{Z}}$ be a function. If  $H(v)$ is allowed for each vertex $v$, then we call $H$ an \emph{allowed function.}
The $H$-factor problem, i.e., determining whether a graph contains $H$-factors, is NP-complete in general.
 For the case when $H$ is an allowed function, \cite{Lovasz} gave a structural description.  In fact, Lov\'asz introduced the definition of negative degree by giving a 2-end-coloring of edges.
By defining the negative degree for a general graph $G$, Lov\'asz may study the degree constrained factor problems of
mixed graphs (including multiple edges, loops, directed edges, two way edges($\leftrightarrow$ or $\rightarrow\leftarrow$ (one edge))).
   \cite{Cor88}
provided the first polynomial time algorithm for $H$-factor problem
with~$H$
being allowed and so give an affirmative answer to the problem proposed by Lov\'asz and Plummer. %

 A classical approach, due to Tutte, for studying $f$-factor problems is to look for
 reductions to the simpler matching problem. For studying $H$-factor problems,
 where every gap of $H(v)$ has the same parity,
\cite{Sz09} used a reduction to local $K_2$ and
factor-critical subgraph packing problem of \cite{CHP}. 
  The idea of reducing a degree prescription to other matching problems  appeared in works of
   \cite{Cor88}. 
    \cite{Cor88} and  \cite{Lo93}  considered reductions to the edge and triangle packing problem,  which can be translated into 1-anti-factor problem.
 Let $G$ be a graph, $U=V(G)$ and let $W$ be the set of all edges and triangles of $G$. Let $G'=(U,W;E')$ be a bipartite graph, where $E'=\{xy\ |\ x\in U,y\in W\mbox{ and }x\in V(y)\}$. Then $G'$ has a 1-anti-factor if and only if $G$ contains a set of vertex-disjoint edges and triangles covering $V(G)$.

  \cite{SH} showed that every graph $G$ contains an $H$-factor  when $|\{1,\ldots,d_G(v)\}-H(v)|=1$ holds for all $v\in V(G)$.  \cite{ADMRT}  showed that every graph $G$ contains a factor $F$ such that $d_F(v)\in \{a_v^-,a_v^-+1,a_v^+,a_v^++1,\}$ for all $v\in V(G)$, where $d_G(v)/3\leq a_v^-\leq d_G(v)/2-1$ and $d_G(v)/2\leq a_v^+\leq 2d_G(v)/3$. \cite{ADR}  slightly improved the result in \cite{ADMRT} and obtained a similar result for  bipartite graphs. For more results on non-consecutive $H$-factor problems of graphs, we refer readers to \cite{LWY,Lu16,TWZ16}.

 However, there is no nice formula to determine whether a bipartite graph  contains a 1-anti-factor. So it is interesting to classify  bipartite graphs with  1-anti-factors.
 Yu and Liu (see \cite{YL09}, Page 76) asked whether every connected $r$-regular bipartite graph contains a 1-anti-factor.
In this paper, we give an affirmative answer to Yu and Liu's problem and obtain the following result.
\begin{theorem}\label{main}
Let $k\geq 3$ be an integer. Every $k$-regular bipartite graph  contains a 1-anti-factor.
\end{theorem}

The rest of the paper is organized as follows. In Section 2, we
introduce Lov\'asz's $H$-Factor Structure Theorem that is needed in
the proof of Theorem 1.3. The proof of Theorem \ref{main} will be presented
in Section 3.

\section{Lov\'asz's $H$-Factor Structure Theorem}\label{sec2}

Let $F$ be a spanning subgraph of $G=(V,E)$ and let $H : V (G)\rightarrow 2^{\mathbb{Z}}$ be an allowed function. Following \cite{Lovasz}, one may measure the ”deviation” of $F$ from the condition (1) by $\boldsymbol{\nabla_H(F,G)}:=\sum_{v\in V(G)} \min\{|d_F(v)-h|: h\in
H(v)\}$. Moreover, let $\boldsymbol{\nabla_H(G)} =\min\{\nabla_H (F, G) : \mbox{ $F$ is a  spanning }$ $\mbox{subgraph of $G$}\}.$ $\nabla_H(G)$ is called \emph{deficiency} of $G$ with respect to the function $H$. The subgraph $F$ is said to be \emph{$H$-optimal} if $\nabla_H(F,G)=\nabla_H (G)$. It is clear that $F$ is an $H$-factor if and only if $\nabla_H(F, G)=0$, and any $H$-factor (if exists) is $H$-optimal. We study $H$-factors of graphs based on Lov\'asz's structural description to the degree prescribed factor problem.

For $v\in V$, we denote by $\boldsymbol{IH(v)}$ the set of degrees of $v$ in all $H$-optimal spanning subgraphs of $G$, i.e., $IH(v):= \{d_F (v) : \mbox{ F is an $H$-optimal spanning subgraph of $G$}\}$. Based on the relation of the sets $IH(v)$ and $H(v)$, one may partition the vertex set $V$ into four classes:
\begin{align*}
  &C_H(G):= \{v\in V : IH(v)\subseteq H(v)\},\\
  &A_H(G):= \{v \in V-C_H (G) : \min IH(v)\geq \max H(v)\},\\
  &B_H(G):=\{v \in V-C_H (G) : \max IH(v)\leq \min H(v)\},\\
  &D_H(G):=V-C_H(G)-A_H(G)-B_H(G).
\end{align*}
When there is no confusion, we omit the reference to $G$. It is clear that the 4-tuple $(A_H, B_H, C_H, D_H )$ is a partition of $V$. A graph $G$ is said to be \emph{$H$-critical} if it is connected and $D_H=V$.
By the definition of $A_H , B_H, C_H$ the following observations hold:
\begin{itemize}
\item[$($*$)$] for every $v\in A_H$, there exists an $H$-optimal graph $F$ such
that $d_F(v)>\max H(v)$;
\item[$($**$)$] for every $v\in B_H$, there exists an $H$-optimal graph $F$ such
that $d_F(v)<\min H(v)$.
\end{itemize}
We will need the following results of \cite{Lovasz}.
\begin{lemma}[\cite{Lovasz}]\label{interval}
Let $G$ be a simple graph and let $H:V(G)\rightarrow 2^{\mathbb{Z}}$ be an allowed function. Let $v\in D_H$.
\begin{itemize}
\item [$($a$)$] $IH(v)$ consists of consecutive integers. 
\item [$($b$)$]  $IH(v)\cap H(v)$ contains no consecutive integers. 
 \end{itemize}
\end{lemma}

Let $R$ be a connected induced subgraph of $G$. Let $H_{R}:V(R)\rightarrow 2^{\mathbb{Z}}$ be a set
function such that $H_{R}(x)=H(x)$ for all $x\in V(R)$.
\begin{lemma}[ \cite{Lovasz}]\label{lov_lem}
Let $G$ be a graph and let $H:V(G)\rightarrow 2^{\mathbb{Z}}$ be an allowed function.
\begin{itemize}
  \item [$($a$)$] $\nabla_H(G)=c(G[D_H])+\sum_{v\in B_H}(\min H(v)-d_{G-A_H}(v))-\sum_{v\in A_H}\max H(v)$. 
  \item [$($b$)$] If $B_H=\emptyset$, then every connected component $R$ of $G[D_H]$ is $H_R$-critical. 
  \item [$($c$)$] $E_G(C_H, D_H) =\emptyset$.
  \item[$($d$)$] If $G$ is $H$-critical, then $\nabla_H(G)=1$.
\end{itemize}
\end{lemma}

\section{The Proof of Theorem \ref{main}}\label{sec3}

\begin{lemma}\label{Lem_ABL}
Let $p\geq 2$ be an integer.
Let $G=(X,Y;E)$ be a bipartite graph. Let $H:V(G)\rightarrow 2^{\mathbb{Z}}$ such that $H(y)=[\max\{d_G(y),p\}]-\{1\}$ for all $y\in Y$  and $H(x)=\{-1,1\}$ for all $x\in X$. Then $A_H\subseteq X$ and $B_H=\emptyset$.
\end{lemma}

 \begin{proof}Firstly, we show that $B_H=\emptyset$. Suppose that $B_H\neq \emptyset$ and let $v\in B_H$. By the definition of $B_H$, if $v\in X$, then $\max IH(v)\leq \min H(v)=-1$, which is impossible. Thus we may assume that $v\in Y$. This implies that $0\leq \max IH(v)\leq \min H(v)=0$. Hence $IH(v)=\{0\}\subseteq H(v)$, which implies
$v\in C_H$, a contradiction.

Next we show that $A_H\subseteq X$ by contradiction. Suppose that there exists a vertex $y\in A_H-X$. Since $p\geq 2$, by the definition of set $A_H$,
  we have that $d_G(y)\geq \max IH(y)\geq\min IH(y)\geq \max H(y)\geq d_G(y)$.
Thus we may infer that $IH(y)=\{d_G(y)\}\subseteq H(y)$, which implies that $y\in C_H$ by the definition, a contradiction. This completes the proof.
\end{proof}

\begin{lemma}\label{H-factor-critical}
Let $p\geq 2$ be an integer.
Let $G=(X,Y;E)$ be a bipartite graph and let $H:V(G)\rightarrow 2^{\mathbb{Z}}$ such that $H(y)=[\max\{d_G(y),p\}]-\{1\}$ for all $y\in Y$  and $H(x)=\{-1,1\}$ for all $x\in X$. If $G$ is $H$-critical, then the following properties hold.
\begin{itemize}
\item [$(i)$] $G-x$ contains an $H_{G-x}$-factor  for all $x\in X$;

\item [$(ii)$] $IH(u)\subseteq \{0,1,2\}$ for  all $u\in V(G)$;

\item [$(iii)$]  $|X|$ is odd;

\item [$(iv)$] Let $y\in Y$ such that $d_G(y)\geq 3$. Then there exist three vertices  $x_1,x_2,x_3\in N_G(y)$ such that $\nabla_{H_{G'}}(G')=2,$ where $G'=G-\{x_1,x_2,x_3,y\}$.
\end{itemize}
\end{lemma}

\begin{proof} Let $G$ be $H$-critical. By the definition of $H$-critical graph
and Lemma \ref{lov_lem} (d), we have that $D_H=V(G)$ and
$\nabla_H(G)=1$. For any $x\in X$, by the definition of $D_H$, there
exists an $H$-optimal subgraph $F$ of $G$ such that $d_F(x)=0$ and
$d_F(w)\in H(w)$ for all $w\in V(G)-\{x\}$. Hence $G-x$ contains an
$H_{G-x}$-factor. This completes the proof of (i).

Next we show (ii). Suppose that there exists a vertex $u\in V(G)$ and an integer $r\geq 3$ with $r\in IH(u)$.
 Since $\nabla_H(G)=1$ and $H(x)=\{-1,1\}$ for any $x\in X$,  we have  $u\in Y$.
From the definition of $D_H$,  we may infer that $IH(u)-H(u)\neq \emptyset$.
Recall that $H(u)=[\max\{d_G(u),p\}]-\{1\}$. Thus we have $1\in IH(u)$.
By Lemma \ref{interval} (a), $IH(u)$ is an interval, which implies  $\{2,3\}\subseteq IH(u)$. Then we have $\{2,3\}\subseteq IH(u)\cap H(u)$, contradicting to Lemma \ref{interval} (b). This completes the proof of (ii).

Given $x\in X$, since $x\in D_H=V(G)$, we may choose an $H$-optimal
subgraph $F$ of $G$ such that $d_F(x)=0$. Note that $\nabla_H(G)=1$.
Thus we have $d_F(w)\in H(w)$ for all $w\in V(G)-\{x\}$.  Since $F$
is bipartite,
\begin{align}\label{Ineq-odd-X}
\sum_{y\in Y}d_F(y)=e_F(X,Y)=\sum_{x\in X}d_F(x)=|X|-1.
\end{align}
By (ii), we have that $d_F(y)\in \{0,2\}$ for all $y\in Y$. So
we have that
 $\sum_{y\in Y}d_F(y)$ is even. By (\ref{Ineq-odd-X}), $|X|$ is odd.  This completes the proof of (iii).

Now we show that (iv) holds. Let $F$ be an $H$-optimal subgraph of
$G$ such that $d_F(y)=1$ and let $N_F(y)=\{x\}$.  Since
$\nabla_H(G)=1$ and $d_F(y)=1\notin H(y)$, we have $d_F(w)\in H(w)$
for all $w\in V(G)-\{y\}$.
Let $x_2,x_3\in N_G(y)-x$. Then we have that $d_{F+x_2y+x_3y}(y)=3\in H(y)$.
One can see that $d_{F+x_2y+x_3y}(w)=d_{F}(w)\in H(w)$ for all $w\in V(G)-\{x_2,x_3,y\}$
and $d_{F+x_2y+x_3y}(x_i)=2$ for $i\in \{2,3\}$.
Set $G'=G-\{y,x,x_2,x_3\}$. Let $y_i\in N_F(x_i)-\{y\}$ for
$i\in\{2,3\}$. (Note that $y_2=y_3$ is possible.) Thus we have
$d_{F-\{y,x,x_2,x_3\}}(w)\in H(w)=H_{G'}(w)$ for all $w\in
V(G')-y_2-y_3$. Recall that $d_{F}(y_i)\in H_{G'}(y_i)$ for $i\in
\{2,3\}$. One can see that
\[
\nabla_{H_{G'}}(F-\{y,x,x_2,x_3\};G')\leq 2.
\]
Hence we have
\[
\nabla_{H_{G'}}(G')\leq 2.
\]
Since $G$ contains no $H$-factors, we have
 \[
\nabla_{H_{G'}}(G')\geq 1.
\]
If  $\nabla_{H_{G'}}(G')= 1$, let $F'$ be an $H_{G'}$-optimal
subgraph of $G'$, then $F'\cup \{xy,x_2y,x_3y\}$ is also an
$H$-optimal subgraph of $G$, which implies $3\in IH(y)$,
contradicting to (ii). This completes the proof. \end{proof}

\begin{theorem}\label{H-factor-iff}
Let $p\geq 2$ be an integer.
Let $G=(X,Y,E)$ be a bipartite graph and let $H:V(G)\rightarrow 2^{\mathbb{Z}}$ such that $H(y)=[\max\{d_G(y),p\}]-\{1\}$ for all $y\in Y$  and $H(x)=\{-1,1\}$ for all $x\in X$.  Then $G$ contains an $H$-factor if and only if for any subset $S\subseteq X$,
we have
\begin{align}\label{q(G-S)>|S|}
q(G-S)\leq |S|,
\end{align}
where $q(G-S)$ denotes the number of connected components $R$ of $G-S$, such
that $R$ is $H_R$-critical.
\end{theorem}

\begin{proof} Firstly, we prove the necessity. Suppose that $G$ contains an $H$-factor $F$.
  Let $R_1,\ldots, R_{q}$ denote these $H_R$-critical components of $G-S$. Since $R_i$ contains no $H_{R_i}$-factors, every $H$-factor  of $G$ contains at least an edge from $R_i$ to $S$. Thus
\[
q(G-S)\leq \sum_{x\in S}d_F(x)= |S|,
\]
which implies $q(G-S)\leq |S|$.

Next, we prove the sufficiency. Suppose that $G$ contains no
$H$-factors. Let $A_H,B_H,C_H,D_H$ be defined as in Section 2. By
Lemma \ref{Lem_ABL}, $A_H\subseteq X$ and $B_H=\emptyset$.

By Lemma \ref{lov_lem} (a), we have
\begin{align*}
0<\nabla_H(G)&=c(G[D_H])+\sum_{v\in
B_H}(\min H(v)-d_{G-A_H}(v))-\sum_{v\in A_H}\max H(v)\\
&=c(G[D_H])-|A_H|,
\end{align*}
i.e.,
\begin{align}\label{DL<AL}
c(G[D_H])>|A_H|.
\end{align}
  By Lemma \ref{Lem_ABL}, we have $B_H=\emptyset$.   By Lemma \ref{lov_lem} (b), every connected component  $R$ of $G[D_H]$ is also $H_R$-critical.
Then, by (\ref{DL<AL}),
\[
q(G-A_H)\geq c(G[D_H])>|A_H|.
\]
This completes the proof. \end{proof}

From the proof of Theorem \ref{H-factor-iff} and Lemma \ref{lov_lem} (b), one can see  the following result.
\begin{lemma}\label{H-anti-factor-if}
Let $p\geq 2$ be an integer.
Let $G=(X,Y,E)$ be a bipartite graph and let $H:V(G)\rightarrow 2^{\mathbb{Z}}$ such that $H(y)=[\max\{d_G(y),p\}]-\{1\}$ for all $y\in Y$  and $H(x)=\{-1,1\}$ for all $x\in X$.  If  $G$ contains no $H$-factors, then
\begin{align}\label{DL<AL}
\nabla_H(G)=c(G[D_H])-|A_H|,
\end{align}
where every connected component $R$ of $G[D_H]$   is $H_R$-critical and also a connected component of $G-A_H$.
\end{lemma}

\begin{lemma}\label{H-factor-H-critical}
Let $k\geq 2$ be an integer. Let $G=(X,Y;E)$ be a connected $k$-regular bipartite graph and let $H:V(G)\rightarrow 2^{\mathbb{Z}}$ such that $H(y)=[k]-\{1\}$ for all $y\in Y$  and $H(x)=\{-1,1\}$ for all $x\in X$. Then either $G$ contains an $H$-factor or $G$ is $H$-critical.
\end{lemma}

\begin{proof}
Suppose that $G$ contains no $H$-factors and is not $H$-critical.
By Lemma \ref{Lem_ABL}, we have that
\begin{align}\label{Ineq:BL}
B_H=\emptyset\ \mbox{and}\ A_H\subseteq X.
\end{align}
Since $G$ is not $H$-critical, we have $D_H\neq V(G)$. Thus we infer that
$A_{H}\neq \emptyset$, otherwise, $C_H=V(G)-D_H\neq \emptyset$ and by Lemma \ref{lov_lem} (c), $E_G(C_H,D_H)=\emptyset$, a contradiction since $G$ is connected.

Recall that $H$ contains no $H$-factors. By Lemmas \ref{Lem_ABL} and
\ref{H-anti-factor-if}, we have $B_H=\emptyset$, $A_H\subseteq X$ and
\begin{align}\label{eq:1}
0<&\nabla_H(G)= c(G[D_H])-|A_H|.
\end{align}
Let $R_1,\ldots, R_q$ denote connected components of $G-A_H$, where   $q=c(G-A_H)$. Since $G$ is a connected regular bipartite graph and $A_H\subseteq X$, we have $|X|=|Y|$ and every connected component $R$ of $G-A_H$ satisfies $|V(R)\cap X|<|V(R)\cap Y|$. So we have
\[
qk\leq k\sum_{i=1}^q(|V(R_i)\cap Y|-|V(R_i)\cap X|)= \sum_{i=1}^q e_G(V(R_i),A_H)= \sum_{x\in A_H}d_G(x)= k|A_H|,
\]
which implies
\[
c(G[D_H])\leq q=c(G-A_H)\leq |A_H|,
\]
contradicting to (\ref{eq:1}). This completes the proof. \end{proof}

Let $\mathcal {H}$ be the set of graphs $G$, which satisfies the
following properties:
\begin{itemize}
\item[$(a)$] $G$ is a connected bipartite graph with color classes $X,Y$;

\item[$(b)$]   $|X|=|Y|-1$;

\item[$(c)$] $d_{G}(x)=3$ for every vertex $x\in X$ and  $d_{G}(y)\leq3$
for every vertex $y\in Y$.
\end{itemize}

\begin{lemma}\label{non-H-critical}
If $G\in \mathcal {H}$, then $G$ is not $H$-critical, where $H:V(G)\rightarrow 2^{\mathbb{Z}}$ is a function such that $H(x)=\{-1,1\}$ for all $x\in X$ and $H(y)=\{0,2,3\}$ for all $y\in Y$.
\end{lemma}

\begin{proof} Suppose that the result does not hold.
 Let $G\in \mathcal {H}$ be an $H$-critical graph with the smallest  order. By Lemma \ref{H-factor-critical} (iii), $|X|$ is odd.
  Recall that $|X|=|Y|-1$ and $d_G(x)=3$ for all $x\in X$. Hence $|X|+1=|Y|\geq 4$ and  there exists $y\in Y$ such that $d_G(y)=3.$
 If
$|Y|=4$, then $|X|=3$ and the spanning subgraph of $G$ with edge  set $\{xy\ |\ x\in N_G(y)\}$ is an $H$-factor, a contradiction.
Hence we may assume that $|X|\geq 5$.

Let $N(y)=\{x_1,x_2,x_3\}$ and $G'=G-N[y]$. Let $H'=H_{G'}$. By Lemma \ref{H-factor-critical} (iv), we have $\nabla_{H'}(G')=2$. Let $A':= A_{H'}(G')$, $B':=B_{H'}(G')$, $C':=C_{H'}(G')$ and $D'=D_{H'}(G')$. By Lemma \ref{Lem_ABL}, $B'=\emptyset$. By Lemma \ref{H-anti-factor-if}, we have
\begin{align}\label{DH'<AH'}
\nabla_{H'}(G')=c(G'[D'])-|A'|= 2.
\end{align}

Now we show that $G'[D']$ contains a connected component $R$ such that $R\in \mathcal{H}$, which contradicts to the choice
of $G$ since $R$ is $H_{R}$-critical and $|V(R)|<|V(G)|$.
Let  $q:=c(G'-A')$.
Let $R_1,\ldots,R_q$ denote the connected components of $G'-A'$.
Note that for every connected component $R$ of $G-A'$, $d_R(x)=3$ for all $x\in V(R)\cap X$. So we have $|V(R)\cap X|<|V(R)\cap Y|.$
Recall that  $|X|=|Y|-1$. Moreover, one can see that $|X|=\sum_{i=1}^q|V(R_i)\cap X|+|A'|+3$ and $|Y|=\sum_{i=1}^q|V(R_i)\cap Y|+1$.
So we may infer that
\begin{align}\label{3.6-eq:2}
\sum_{i=1}^q|V(R_i)\cap X|+|A'|+3=\sum_{i=1}^q|V(R_i)\cap Y|\geq \sum_{i=1}^q|V(R_i)\cap X|+q,
\end{align}
i.e.,
\begin{align}\label{q>|A'|+3}
q\leq |A'|+3.
\end{align}
Since $E_{G'}(C',D')=\emptyset$, combining (\ref{DH'<AH'}), we have $q\geq c(G[D'])= |A'|+2\geq 2$. So $q\in \{|A'|+2,|A'|+3\}$.
By (\ref{3.6-eq:2}), each connected component $R$ of $G'-A'$ except at most one satisfies $|V(R)\cap X|=|V(R)\cap Y|-1$. Since $c(G[D'])\geq 2$, we have $G[D']$ contains an $H_R$-critical component $R$ such that $|V(R)\cap X|=|V(R)\cap Y|-1$. By Lemma \ref{H-factor-critical} (iii), $|V(R)\cap X|$ is odd and so $V(R)\cap X\neq \emptyset$. Hence we have $R\in \mathcal {H}$. This completes the proof. \end{proof}

\noindent \textbf{Proof of Theorem \ref{main}:}
  Let $G$ be a $k$-regular bipartite graph with bipartition $(X,Y)$.
 Let $H:V(G)\rightarrow 2^{\mathbb{Z}}$ such that $H(x)=\{-1,1\}$ for all $x\in X$ and $H(y)=\{0,2,3\}$ for all $y\in Y$.
  Clearly, if $G$ has an $H$-factor, then $G$ has a 1-anti-factor. By Hall's Theorem,
$G$ contains a $3$-factor. Thus it is sufficient for us to show that
every connected 3-regular bipartite graph contains an $H$-factor. So
we may assume that $G$ is a connected 3-regular bipartite graph.
 By contradiction, suppose that $G$ contains no $H$-factors.

 By Lemma \ref{H-factor-H-critical}, we may assume that
$G$ is $H$-critical. 
Let $y\in Y$ and $G'=G-N[y]$.  Let $H':=H_{G'}$, $D':=D_{H'}(G')$, $A':=A_{H'}(G')$, $B':=B_{H'}(G')$  and $C':=C_{H'}(G')$.
By Lemma \ref{H-factor-critical} (ii) and (iv), we have  that $IH(y)\subseteq \{0,1,2\}$ and  $\nabla_{H'}(G')=2$.
By Lemma \ref{H-anti-factor-if}, we have
\begin{align}\label{th-eq:1}
2=\nabla_{H'}(G')
=c(G'[D'])-|A'|,
\end{align}
By Lemma \ref{Lem_ABL}, we have $B'=\emptyset$.
  Let $q:=c(G'-A')$. Let $R_1, \ldots, R_q$ be the connected components of $G'-A'$.

Now we will show that $G'[D']$ contains a connected component $R$ such that $R$ is $H_R$-critical and $R\in \mathcal{H}$, which contradicts to Lemma \ref{non-H-critical}.  (The proof is completely similar with that of Lemma \ref{non-H-critical}.)
Note that $|X|=|Y|$, $|X|=\sum_{i=1}^q |V(R_i)\cap X|+3+|A'|$ and
\begin{align}\label{th-eq:2}
|Y|=\sum_{i=1}^q |V(R_i)\cap Y|+1\geq \sum_{i=1}^q |V(R_i)\cap X|+q+1.
\end{align}
So we have $q\leq |A'|+2$. By (\ref{th-eq:1}), we have $q\geq c(G'[D'])= |A'|+2$. Thus $q=|A'|+2$ and so the equality  holds for (\ref{th-eq:2}), which implies that for every connected component $R$ of $G'-A'$, it is $H_R$-critical and $|V(R)\cap X|=|V(R)\cap Y|-1$.  So every connected component of $G'[D']$ belongs to $\mathcal{H}$.
This completes the proof of Theorem 1.2. \qed

\noindent \textbf{Remark 1.}
The bound that $k\geq 3$ in Theorem \ref{main} is sharp. Let $m\in \mathrm{N}$ be a positive integer. For
example, $C_{4m+2}$ is a 2-regular graph and contains no $H$-factors. However, it is easy to show that $C_{4m}$ contains
an $H$-factor.

\noindent \textbf{Remark 2.} Theorem  \ref{main} does not hold for multi-graphs. By doubling every second edge in $C_{4m+2}$, we get a 3-regular bipartite multi-graph $G$. But, as one sees in Remark 1 that $C_{4m+2}$ does not contain an $H$-factor, one sees that neither does $G$.

\acknowledgements
\label{sec:ack}
The authors would like to thank the anonymous Reviewers for all valuable comments and suggestions to greatly improve the quality of our paper.

\nocite{*}
\bibliographystyle{abbrvnat}
\bibliography{dmtcs-6463}
\label{sec:biblio}

\end{document}